\documentstyle[amsfonts]{article}
\newcommand{\E}{{\Bbb E}}

\newcommand{\R}{\Bbb R}

\newtheorem{theorem}{Theorem}[section]

\newtheorem{corollary}[theorem]{Corollary}

\begin{document}

\title{Bismut type formulae for differential forms}
\author{K.D. Elworthy${}^*$\and Xue-Mei Li${}^{**}$}
\date{}
\maketitle

\footnote{
 *Mathematics Institute, University of Warwick, Coventry CV4 7AL, UK.
Research supported by EC programme Science plan ERB 4002PL910459.

\hskip 8pt ** Department of mathematics, University of Connecticut,
Storrs, USA  and Mathematical Sciences Research Institute, Berkeley, 
California.  Research supported by NSF
 grant DMS-9626142 and Alexander Humboldt-Stiftung.}

\begin{abstract}
   Formulae are given for $dP_t \phi$, $d^*P_t\phi$ and 
 $\Delta P_t\phi$ for $P_t$ the heat semigroup acting on a q-form $\phi$.
 The formulae are Brownian motion expectations  of $\phi$ composed with
 random  translations determined by Weitzenbock curvarure terms.
Derivatives of the curvature are not involved.
\end{abstract}
{\large Formules de Bismut sur les formes diff\'erentielles}

Resum\'e:  Nous donnons les expressions de   $dP_t \phi$, $d^*P_t\phi$ et 
 $\Delta P_t\phi$, o\`u $P_t$ d\'esigne le semi-groupe  de  la chaleur
 agissant sur les  formes diff\'erentielles d'ordre $q$. Elles s'\'ecrivent
comme l'esp\'erance de la fonctionnelle  $\phi$ compos\'ee avec les
 transports  parall\`eles al\'eatoires associ\'es aux  courbures de 
Weitzenbock. Les  d\'eriv\'ees de la courbure n'interviennent pas dans ces
 formules.
\bigskip

{\Large Version francaise abr\'eg\'ee}

\noindent
{\bf A.} Soit $M$  une vari\'et\'e compacte, $\Delta$  le laplacien
 avec la convention de signe $\Delta=-(d\delta+\delta d)$ et
 ${\cal  R}:\wedge^*TM \to \wedge^*TM$ est la courbure de  Weitzenbock.
 On a la formule de Weitzenbock
\begin{equation}  \label{f-1}
\Delta \phi=\hbox{trace} \nabla^2\phi-\phi\circ {\cal R}
\end{equation}
On note  $ {\cal R}^q$, etc. les restrictions des op\'erateurs sur
 les formes diff\'erentielles d'ordre $q$ et sur les $q$-vecteurs.
 Rappelons ${\cal R}^1=Ric^{\#}$, le tenseur de Ricci.
  On d\'esigne par $\{P_t: t\ge 0\}$ le semi-groupe  
 ${\ {\rm e}^{{\ {\frac{1}{2}} } t \Delta} }$ agissant sur les formes
 diff\'erentielles de $L^2$ (par   auto-adjonction essentielle).
 Il agit donc sur les formes continues par 
la formule de Feyman-Kac d'Airault $P_t\phi=\E\phi(V_t),$ 
o\`u $V_t$ est la solution de l'\'equation   covariante
 le long des trajectoire  browniennes $\{x_t:t\ge 0\}$ 
\begin{equation}
{\frac{D V_t}{\partial t}}=-{\frac 12}{\cal R}(V_t),\hskip 18ptV_0\in
\wedge _{x_0}^{*}TM.  \label{f-3}
\end{equation}
On d\'efinit $W_t:\wedge^*TM\to \wedge^*TM$ par   $V_t=W_t(V_0)$.
 Soit  $\{\breve{B}_t: t\ge 0\}$ l'anti-d\'eveloppement
 stochastique du  mouvement brownien $\{x_t:t\ge 0\}$, et $//_t$ le
transport parall\`ele  de la connexion de Levi-Civita le long 
des trajectoires.


La formule de Bismut est pour les fonctions born\'ees  mesurables
 $f:M\to \R$:
\begin{equation}\label{f-4}
dP_tf(v_0)=\frac 1t{\Bbb E\;}f(x_t)\int_0^t<//_sd\breve{B}_s,W_s>.
\end{equation}
On a donc sa  formule pour le  gradient du noyau  de la chaleur 
 $p_t(x,y)$ sur les fonctions\cite{Bismut84}:
\begin{equation}
\nabla _x\log p_t(x_0,y)={\Bbb E} \left\{\frac 1t
\int_0^tW_s^{*}(//_s d\breve{B}_s) \left | x_t=y \right.\right\}.
\label{f-5}
\end{equation}

Nous g\'en\'eralisons ces formules aux cas des $q$-formes, $1\le q< n$.

\noindent
{Th\'eor\`eme 1.}
{\it  Soit $\phi$  une forme diff\'erentielle, born\'ee 
 mesurable sur la vari\'et\'e $M$. 
\begin{equation}
d(P_t\phi )\;(V_0)=\frac 1t{\Bbb E\;\phi }\left(
W_t^q\int_0^t(W_s^q)^{-1}\,\iota _{//_sd\breve{B}_s}W_s^{q+1}(V_0)\right),
\hskip 16pt  V_0\in \wedge^*T_{x_0}M,
\label{f-6}
\end{equation}
o\`u  $\iota _{//_sd\breve{B}_s}$ d\'esigne le produit int\'erieur
 (l'op\'erateur d'annihilation):
\[
\iota _{//_sd\breve{B}_s}v^1\wedge \dots \wedge
v^{q+1}=\sum_{j=1}^{q+1}(-1)^{j+1}<//_sd\breve{B}_s,v^j>v^1\wedge\dots\wedge
\widehat{v^j}\dots v^{q+1}.
\]  }

Nous avons un corollaire analogue \`a (\ref{f-5}) :

{Corollaire 2.} {\it
Pour le noyau de la chaleur $p_t^q (x_0,y)$ d\'efini sur les $q$-formes:
\begin{equation}
d_xp_t^q(x_0,y)(V_0)=\frac 1t p_t^q(x_0,y) \left(
E\, \{\,W_t^q\int_0^t(W_s^q)^{-1}\,\iota _{//_sd\breve{B}_s}W_s^{q+1}(V_0) 
\left| x_t=y_0 \right.\, \}\,\right). 
\label{f-7}
\end{equation}}
Nous remarquons que ces formules ne font pas intervenir les d\'eriv\'ees
du  tenseur de courbure. Dans  \cite{Norris93} Norris donne une 
 expression des  d\'eriv\'ees covariantes de tel semi-groupes
 de la chaleur.
Toutefois ses formules font intervenir les d\'eriv\'ees du tenseur
de courbure.

La formule (\ref{f-6}) provient d'une  formule non-intrins\`eque de
 Li \cite{thesis}:
\begin{equation}
d(P_t^q\phi )=\frac 1tE\,\int_0^t<//_{s\,}d\breve{B}_s,T\xi _s(-)>\wedge
\phi \left( \wedge ^qT\xi _t(-)\right) ,  \label{f-8}
\end{equation}
voir aussi  Elworthy-Li \cite{EL-LI}, en utilisant les techniques de
 Elworthy-Yor \cite{Elworthy-Yor}, comme dans
 Elworthy-Li \cite{EL-LI.ibp}. Pour les details voir la version anglaise.

On peut d\'eduire le Th\'eor\`eme 1 d'une formule 
d'int\'egration par partie  plus g\'en\'erale sur les espaces de
 chemins pour les $q$-formes, voir \cite{EL-LI98}. Quand on
remplace le mouvement brownien par le mouvement brownien avec 
d\'erive la m\^eme d\'emonstration fournit des formules avec des changements
\'evidents \cite{EL-LI}, \cite{EL-LI.ibp}. De m\^eme pour le changement 
 de Thalmaier, o\'u
 ${1\over t}\int_0^t$ est remplac\'e par 
${1\over h}\int_{\delta}^{\delta+h}$, 
 $0\le\delta <\delta+h\le t$  \cite{EL-LI}.

\bigskip

\noindent
{\bf B.} \underline{Les formules pour $d^*P_t\phi$.} 
Nous donnons un autre corollaire, suite \`a une conversation avec Bruce 
Driver. Driver   a d\'emontr\'e la cas particulier de $q=1$ par une
 m\'ethode diff\'erente \cite{Driver}. Il a
 \'egalement prouv\'e
  une formule intriguante faisant intervenir 
 une int\'egrale stochastique backward \cite{Driver97}.

\noindent{Corollaire 3.} {\it
Soit  $\phi$ une forme diff\'erentielle d'ordre $q$ de $L^2$. On a
\begin{equation}\label{f-divergence}
d^*P_t\phi= -{1\over t}
 \E\,\phi \left( W_t^q\int_0^t (W_s^q)^{-1}
\left(//_s d\breve B_s\wedge W_s^{q-1} (-)\right)\right).
\end{equation}}
 Voir la version anglaise pour une preuve.
\bigskip

\noindent
{\bf C.} \underline{Les formules pour $\Delta P_t\phi$.}
La methode utilis\'ee dans \cite{EL-LI} pour obtenir les d\'eriv\'ees
d'ordre sup\'erieur du semi-groupe de  la chaleur agissant sur les fonctions
peut \^etre utilis\'ee ici pour obtenir des r\'esultats plus
 \'el\'egants:

{Corollaire 4.} {\it
\begin{equation}
\Delta P_t\phi (V)={4\over t^2}\E\,
\phi\left(  W_t \int_{t/2}^t(W_s)^{-1}dM_s(V)\right), \hskip
36pt \hbox{o\`u}
\end{equation}
\begin{eqnarray*}
dM_s(V)&=& -
\iota_{//_s d \breve B_s} W_s\left( \int_0^{t/2} (W_r)^{-1}
\left( //_r d\breve B_r \wedge W_r(V)\right)\right)\\
&&   -
//_s d \breve B_s \wedge W_s \left(\int_0^{t/2} (W_r)^{-1} 
\left(\iota_{//_r d\breve B_r} W_r(V)\right)\right).
\end{eqnarray*}}

Nous pouvons it\'erer cette m\'ethode afin d'obtenir les formules pour
  $(\Delta)^r P_t\phi, r=1,2\dots $. Elles permettent d\'etablir
les formules pour $d_x p_t^q(\cdot,y)$, $d^*_x p_t^q(\cdot, y)$, 
$(\Delta_x)^r p_t^q(\cdot, y)$. En particulier on en deduit que 
$d_x p_t^q(\cdot,y)$, $d^*_x p_t^q(\cdot, y)$, $(\Delta_x)^r p_t^q(\cdot, y)$
et  $dP_t\phi$, $d^*P_t\phi$, $(\Delta)^r P_t \phi$, $r=1,2,\dots$ sont
uniformement born\'es par une constante
 d\'ependant uniquement de $t$, des courbures de Weitzenbock  ${\cal R}^q$, 
${\cal R}^{q\pm 1}$ et du maximum  de $|\phi|$.

\pagebreak

{\bf A.}
Let $M$ be a compact Riemannian manifold with Hodge-deRham Laplacian $%
\Delta $ given the sign convention $\Delta=-(d\delta+\delta d)$. There is
then the Weitzenbock formula 
\begin{equation}  \label{1}
\Delta \phi=\hbox{trace} \nabla^2\phi-\phi\circ {\cal R}
\end{equation}
where ${\cal R}:\wedge^*TM \to \wedge^*TM$ is the Weitzenbock curvature
 term. Let ${\cal R}^q$, etc. denote the restriction to 
$q$-forms and $q$-vectors. Recall ${\cal R}^1=Ric^{\#}$, the Ricci tensor.
 Let $\{P_t:t\ge 0\}$ be the semigroup ${\ {\rm e}^{{\ {\frac{1}{2}} } t
\Delta} }$ acting on $L^2$ forms (by the essentially self-adjointness of $%
\Delta$) and on bounded measurable forms via Airault's version of the
 Feyman-Kac formula
\begin{equation}
P_t\phi (V_0)={\Bbb E}\phi (V_t)  \label{2}
\end{equation}
where $V_t$ satisfies the covariant equation along Brownian paths $%
\{x_t:t\ge 0\}$ 
\begin{equation}
{\frac{D V_t}{\partial t}}=-{\frac 12}{\cal R}(V_t),\hskip 18ptV_0\in
\wedge^* T_{x_0}TM.  \label{3}
\end{equation}
We will write $V_t=W_t(V_0)$ to give $W_t: \wedge^*TM\to \wedge^* TM$. 
Let $\{\breve{B}_t: t\ge 0\}$ be the stochastic
anti-development of $\{x_t:t\ge 0\},$ a Brownian motion
on $T_{x_0}M$, and denote by $//_t$ the parallel translation for the
Levi-Civita connection along our  paths.

\bigskip

Bismut's formula is for bounded measurable functions $f:M\to {\Bbb R}$.
It is
\begin{equation}\label{4}
dP_tf(v_0)=\frac 1t{\Bbb E\;}f(x_t)\int_0^t<//_sd\breve{B}_s,W_s>
\end{equation}
which leads to his gradient formula for the heat kernel $p_t(x,y)$ on
functions \cite{Bismut84}: 
\begin{equation}
\nabla _x\log p_t(x_0,y)={\Bbb E} \left\{\frac 1t
\int_0^tW_s^{*}(//_s d\breve{B}_s) \left | x_t=y \right.\right\}.
\label{5}
\end{equation}

Here we extend these formulae to the case of $q$-forms, $1\le q\le n$ to show:
\begin{theorem}
\label{th:1}
For $\phi $ a bounded measurable $q$-form on $M$ and
 $V_0\in \wedge^*T_{x_0}M$
\begin{equation}
d(P_t\phi )\;(V_0)=\frac 1t{\Bbb E\;\phi }\left(
W_t^q\int_0^t(W_s^q)^{-1}\,\iota _{//_sd\breve{B}_s}W_s^{q+1}(V_0)\right) ,
\label{6}
\end{equation}
where $\iota _{//_sd\breve{B}_s}$ denotes the interior product (annihilation
operator) given by 
\[
\iota _{//_sd\breve{B}_s}v^1\wedge \dots \wedge
v^{q+1}=\sum_{j=1}^{q+1}(-1)^{j+1}<//_sd\breve{B}_s,v^j>v^1\wedge\dots\wedge
\widehat{v^j}\dots v^{q+1}.
\]
\end{theorem}
As for (\ref{5}) there is the corollary:
\begin{corollary}
For the heat kernel $p_t^q(x_0,y)$ on $q$-forms
\begin{equation}
d_xp_t^q(x_0,y)(V_0)=\frac 1t p_t^q(x_0,y) \left(
E\, \{\,W_t^q\int_0^t(W_s^q)^{-1}\,\iota _{//_sd\breve{B}_s}W_s^{q+1}(V_0) 
\left| x_t=y_0 \right.\, \}\,\right). 
\label{7}
\end{equation}
\end{corollary}
Note these formulae do not involve derivatives of the curvature tensor.
Expressions for covariant derivatives of such heat semigroups have been
obtained by Norris \cite{Norris93}, however these do involve such 
derivatives.

\bigskip

Formula (\ref{6}) is derived from a non-intrinsic formula in Li 
\cite{thesis}:
\begin{equation}
d(P_t^q\phi )=\frac 1tE\,\int_0^t<//_{s\,}d\breve{B}_s,T\xi _s(-)>\wedge
\phi \left( \wedge ^qT\xi _t(-)\right) ,  \label{8}
\end{equation}
see also Elworthy-Li \cite{EL-LI}, by techniques of Elworthy-Yor 
\cite{Elworthy-Yor}, as used in Elworthy-Li \cite{EL-LI.ibp}. The details
 follow:

{\bf Proof of Theorem \ref{th:1}:}
Take an isometric embedding $\alpha :M\rightarrow R^m$, some $m$, using
Nash's Theorem. Let $X(x):R^m\rightarrow T_xN$ be the orthogonal projection,
identifying $T_xM$ with its image in $R^m$ under $d\alpha $. Let $%
\{B_t:t\ge 0\} $ be a Brownian motion on $R^m$. The {\it gradient Brownian
system}
\begin{equation}
dx_t=X(x_t)\circ dB_t  \label{9}
\end{equation}
has a solution flow $\{\xi _t(x):t\ge 0,x\in M\}$ consisting of random
diffeomorphisms $\xi _t:M\rightarrow M$ with derivative denoted by $%
T_{x_0}\xi _t:T_{x_0}M\rightarrow T_{x_t}M$ for $x_t=\xi _t(x_0)$. Each
solution $(x_t:t\ge 0)$ is a Brownian motion on $M$.
Define $J_t:\wedge ^{q+1}TM\to \wedge ^qTM$ by:
\begin{equation}
J_t(V)=\frac 1t\int_0^t<T\xi _r-,//_rd\breve{B}_r>\wedge (\Lambda
^qT\xi _t(-)), \hskip 18pt t>0.  \label{forms:3}
\end{equation}
We see that formula (\ref{8}) can be written as: 
\[
d(P_t\phi )(V)=E\phi (J_t(V)), \hskip 36pt t>0
\]
and $J_t$ satisfies
\[
DJ_t=\iota _{X(x_t)dB_t}\Lambda ^{q+1}T\xi _t+d\Lambda ^q\left( \nabla
X(-)dB_t\right) -{\frac 12}{\cal R}^q(J_t)dt
\]
where the term involving ${\cal R}^q$ comes from the Stratonovich 
correction. See Elworthy-Yor \cite{Elworthy-Yor}. Set
\[
\bar{J}_t=E\left\{ J_t|{\cal F}_t^{x_0}\right\} ,
\]
where ${\cal F}_\cdot^{x_0}$ is the filtration generated by $\{x_\cdot\}$
 and the conditional expectation is understood to be
 $//_t E\left\{ //_t^{-1}J_t\,|\,%
{\cal F}_t^{x_0}\right\} $. Using the decomposition of $B_t$ into relevant
and 'redundant' noise as in Elworthy-Yor \cite{Elworthy-Yor} we have 
\[
d\Lambda ^q\nabla X(-)dB_t=d\Lambda ^q\nabla X(-)/\tilde{/}_td\beta _t
\]
where $/\tilde{/}_t,t\ge 0$ consists of random isometries of $R^m$ and $%
\left\{ \beta _t:t\ge 0\right\} $ is a Brownian motion on
 $\displaystyle{\hbox{Ker}X(x_0)}$
independent of ${\cal F}_t^{x_0}$. From this $\bar{J}_t$ satisfies
\begin{equation}
D\bar{J}_t=\iota _{X(x_t)dB_t}W_t^{q+1}-{\frac 12}{\cal R}^q(\bar{J}_t)dt
 \label{12}
\end{equation}
For a general discussion see Elworthy-LeJan-Li \cite{EL-LJ-LI98}. Solve the
 equation (\ref{12}) to get:
\[
\bar{J}_t=W_t^q\int_0^t(W_s^q)^{-1}\left( \iota _{X(x_s)dB_s}\right)
W_s^{q+1}.
\]
The required result follows. \hfill\rule{3mm}{3mm}

  Theorem \ref{th:1} can also be deduced from a more general integration
 by parts formula for $q$-forms on path space, see \cite{EL-LI98}.
The obvious variation when the Brownian motion is replaced by one with a drift 
$Z$ also holds with the same proof: here $\{W_\cdot\}$ is replaced by 
$\{W_\cdot^Z\}$  as in \cite{EL-LI}, \cite{EL-LI.ibp}.
So also does Thalmaier's variant when 
$ {1\over t}\int_0^t$  is replaced by 
${1\over h}\int_{\delta}^{\delta+h}$ for
 $0\le\delta <\delta+h\le t$ in
 (\ref{8}), (\ref{6}), and (\ref{7}).

\bigskip
 
\noindent
{\bf B.} \underline{Formulae for $d^*P_t\phi$.} 
We give the next corollary following conversations with Bruce Driver. Driver
proved it in the special case $q=1$, by a different method \cite{Driver}.
Also in Driver \cite{Driver97} he gave an intriguing different formula 
in that case involving a backward stochastic integral.
\begin{corollary}
Let $\phi$ be a $L^2$ $q$-form. Then
\begin{equation}\label{divergence}
d^*P_t\phi= - {1\over t}
 \E\,\phi \left( W_t^q\int_0^t (W_s^q)^{-1}
\left(//_s d\breve B_s\wedge W_s^{q-1} (-)\right)\right).
\end{equation}
\end{corollary}
{\bf Proof.}
We can assume that $\phi$ is smooth. Also by going to the double cover if
 necessary we can assume that $M$ is oriented with volume element $d vol$
considered as a $n$-form. The Hodge star construction then gives
$*:\wedge T^*M \to \wedge T^*M$ and
$* :\wedge TM \to \wedge TM$ such that
$\displaystyle{  \alpha\wedge *\beta\,  =<\alpha,\beta> d \,vol}$ and
    $\displaystyle{ v_1\wedge *  v_2=<v_1, v_2> (d\,{vol})^{\#} }$.
Then  $\displaystyle{*\phi(V)=(-1)^{q(n-q)}\phi(*  V)}$ 
for   $V\in \wedge^q T_xM$,   $\displaystyle{ **=(-1)^{q(n-q)}Id}$,  and 
$\displaystyle{d^*\phi =(-1)^{n(q-1)+1}*d*\phi}$.
Furthermore $ *P_t=P_t*$, $\displaystyle{*W_t=W_t * }$,
and $\displaystyle{ *  \iota_u V=(-1)^{q-1} u\wedge (*  V)}$,
 any  $u\in T_xM$.    Thus
$$d^*P_t\phi(-)=(-1)^{n(q-1)+1}*d*(P_t\phi)(-)= (-1)^{q} d*(P_t \phi)(* -)
= (-1)^{q} d(P_t *\phi)(* -).$$
Consequently  (\ref{divergence}) follows from
\begin{eqnarray*}
d^*P_t\phi(-) &=& (-1)^{q}   {1\over t}  \E (*\phi) \left(W_t^{n-q}\int_0^t
  (W_s^{n-q})^{-1} \iota_{X(x_s)dB_s }W_s^{n-q+1}(* -)\right)\\
&=&(-1)^q {1\over t} (-1)^{q(n-q)}\E\, \phi\left( W_t^q\int_0^t (W_s^q)^{-1}
 \left[(-1)^{n-q} X(x_s)dB_s \wedge W_s^{q-1}(* * -)\right] 
\right)\\
&=& -{1\over t} \E\, \phi \left(  W_t^q\int_0^t (W_s^q)^{-1}
\left( X(x_s)dB_s \wedge W_s^{q-1}(-)\right) \right).\\
\end{eqnarray*}
 \hfill\rule{3mm}{3mm}

\noindent
{\bf C.}  \underline{Formulae for $\Delta P_t\phi$.}
The method used in \cite{EL-LI} to obtain higher order derivatives of heat
semigroups on functions can be applied to the present situation, with neater
results. First note that for a $(q+1)$-form $\psi$ equation
 (\ref{6}) at time ${t/2}$ is equivalent to
\begin{equation}\label{14}
dP_{t/2}\psi(V)={2\over t} \E\, \psi_{\xi_t^{t/2}(x_0)}
\left( W_t^{t\over 2} \int_{t/2}^t (W_s^{t/2})^{-1}
\iota_{X(\xi_s^{t/2}(x_0))dB_s} W_s^{t/2}(V)\right)
\end{equation}
where $\xi_s^{t/2}$ and $W_s^{t/2}$, $t/2\le s<\infty$ refer to the
 respective flows starting at time $t/2$ and evaluated at time $s$. Since 
$$d^*dP_t\phi(V)=dP^*_{t/2}(dP_{t/2})\phi(V)$$
equation (\ref{14}), (\ref{8}), the Markov property and the semigroup 
 property of stochastic flows and the development map yield
$$    d^*dP_t\phi(V) =-{4\over t^2}\,
\E \,\phi_{x_t}\left(W_t \int_{t/2}^t (W_s)^{-1}
\iota_{//_s d\breve B_s} W_s
\int_0^{t/2}(W_r)^{-1} //_r d\breve B_r\wedge W_r( V)\right).
$$
Adding on the corresponding formula for $d^*d P_t\phi(V)$ we obtain
\begin{corollary}
\begin{equation}
\Delta P_t\phi (V)={4\over t^2}\E\,
\phi\left(  W_t \int_{t/2}^t(W_s)^{-1}dM_s(V)\right),
 \hskip 36pt \hbox{where}
\end{equation}
\begin{eqnarray*}
dM_s(V)&=& -
\iota_{//_s d \breve B_s} W_s\left( \int_0^{t/2} (W_r)^{-1}
\left( //_r d\breve B_r \wedge W_r(V)\right)\right)\\
&&-
//_s d \breve B_s \wedge W_s \left(\int_0^{t/2} (W_r)^{-1} 
\left(\iota_{//_r d\breve B_r} W_r(V)\right)\right).
\end{eqnarray*}
\end{corollary}
This can be iterated to give formulas for $(\Delta)^r P_t\phi, r=1,2\dots $.
These together with the previous  ones lead to formulae for the corresponding
operators acting on the heat kernels, as for (\ref{7}), and show that
$d_x p_t^q(\cdot,y)$, $d^*_x p_t^q(\cdot, y)$, $(\Delta_x)^r p_t^q(\cdot, y)$
and $dP_t\phi$, $d^*P_t\phi$, $(\Delta)^r P_t \phi$, $r=1,2,\dots$ all have
uniform bounds for fixed positive $t$, depending only on $t$ and bounds
for the Weitzenbock curvatures ${\cal R}^q$, ${\cal R}^{q\pm 1}$, (and
the uniform bound of $\phi$).


\begin{thebibliography}{10}

\bibitem{Bismut84}
J.~M. Bismut.
\newblock {\em Large deviations and the {M}alliavin calculus. {P}rogress in
  {M}ath. 45.}
\newblock Birkh{\"a} user, 1984.

\bibitem{Driver97}
B.~K. Driver.
\newblock Integration by parts for heat kernel measures revisited.
\newblock {\em J. Math. pures et applique\'ee}, 76:703--73, 1997.

\bibitem{Driver}
B.K. Driver.
\newblock Manuscript, 1996.

\bibitem{EL-LJ-LI98}
K.~D. Elworthy, Yves LeJan, and Xue-Mei Li.
\newblock On the geometry of diffusion operators and stochastic flows.
\newblock MSRI prepublication No. 1998-31., 1998.

\bibitem{Elworthy-Yor}
K.~D. Elworthy and M.~Yor.
\newblock Conditional expectations for derivatives of certain stochastic flows.
\newblock In J.~Az\'ema, P.A. Meyer, and M.~Yor, editors, {\em Sem. de Prob.
  XXVII. Lecture Notes in Maths. 1557}, pages 159--172. Springer-Verlag, 1993.

\bibitem{EL-LI98}
K.D. Elworthy and Xue-Mei Li.
\newblock Integration by parts formulae for forms on path spaces of riemannian
  manifolds.
\newblock In preparation.

\bibitem{EL-LI}
K.D. Elworthy and Xue-Mei Li.
\newblock Formulae for the derivatives of heat semigroups.
\newblock {\em J. Funct. Anal.}, 125(1):252--286, 1994.

\bibitem{EL-LI.ibp}
K.D. Elworthy and Xue-Mei Li.
\newblock A class of integration by parts formulae in stochastic analysis {I}.
\newblock In {\em {I}t\^o's Stochastic Calculus and Probability Theory
  (dedicated to Prof. It\^o on the occasion of his eightieth birthday)}.
  Springer-Verlag, 1996.

\bibitem{thesis}
Xue-Mei Li.
\newblock Stochastic flows on noncompact manifolds.
\newblock Ph.D. thesis, University of Warwick, 1992.

\bibitem{Norris93}
J.~Norris.
\newblock Path integral formulae for heat kernels and their derivatives.
\newblock {\em Probability Theory and Related Fields}, 94:525--541, 1993.

\end{thebibliography}
\end{document}